\documentclass[a4paper,10pt]{amsart}
\usepackage{amssymb}

\def\bfit{\bfseries\itshape}

\input xypic
\xyoption{all}
\xyoption{arc}

\newtheorem{theo}{Theorem}[section]

\newtheorem{prop}[theo]{Proposition}

\newtheorem{lem}[theo]{Lemma}

\newtheorem{coro}[theo]{Corollary}

\def\remark#1{{\refstepcounter{theo}\label{#1}\noindent\sc Remark  
\arabic{section}.\arabic{theo} - }}
\def\example#1{{\refstepcounter{theo}\label{#1}\noindent\sc Example 
\arabic{section}.\arabic{theo} - }}
\def\examples#1{{\refstepcounter{theo}\label{#1}\noindent\sc Examples 
\arabic{section}.\arabic{theo} - }}

\def\equat{\refstepcounter{theo}$$~}
\def\endequat{\leqno{\boldsymbol{(\arabic{section}.\arabic{theo})}}~$$}

\newcounter{soussection}[section]

\def\soussection#1{\refstepcounter{soussection}
\noindent{\bfit{\arabic{section}.\Alph{soussection}. #1.}}}

    \def\CM{{\mathbb{C}}}

    \def\NM{{\mathbb{N}}}

    \def\QM{{\mathbb{Q}}}
    
\def\SG{{\mathfrak S}}

    \def\ZM{{\mathbb{Z}}}

    \def\CC{{\mathcal{C}}}

    \def\FC{{\mathcal{F}}}

    \def\IC{{\mathcal{I}}}

    \def\PC{{\mathcal{P}}}
    
    \def\RC{{\mathcal{R}}}

  \def\arm{{\mathrm{a}}}  
  \def\brm{{\mathrm{b}}}

\def\a{\alpha}
\def\b{\beta}
\def\g{\gamma}

\def\D{\Delta}

\def\l{\lambda}
\def\L{\Lambda}

\def\O{\Omega}

\def\s{\sigma}

\def\th{\theta}

\def\t{\tau}

\def\lamb{{\boldsymbol{\lambda}}}       
\def\Lamb{{\boldsymbol{\Lambda}}}

      \def\varpit{{\tilde{\varpi}}}

\DeclareMathOperator{\Hom}{{\mathrm{Hom}}}

\DeclareMathOperator{\Id}{{\mathrm{Id}}}

\DeclareMathOperator{\im}{{\mathrm{Im}}}

\DeclareMathOperator{\Ind}{{\mathrm{Ind}}}

\DeclareMathOperator{\Irr}{{\mathrm{Irr}}}
\DeclareMathOperator{\Ker}{{\mathrm{Ker}}}

\DeclareMathOperator{\Rad}{{\mathrm{Rad}}}

\DeclareMathOperator{\Res}{{\mathrm{Res}}}

\def\to{\rightarrow}
\def\longto{\longrightarrow}

\def\fonction#1#2#3#4#5{\begin{array}{rccc}
{#1} : & {#2} & \longto & {#3} \\
& {#4} & \longmapsto & {#5} 
\end{array}}

\def\vide{\varnothing}

\def\DS{\displaystyle}
\def\SS{\scriptstyle}

\def\finl{~$\SS \square$}

\def\infspe{\hspace{0.1em}\mathop{\preccurlyeq}\nolimits\hspace{0.1em}}

\def\matrice#1{\left(\begin{array}{ccccccccccccccccccc}#1\end{array}\right)}

\def\lexp#1#2{\kern\scriptspace\vphantom{#2}^{#1}\kern-\scriptspace#2}
\def\le{\hspace{0.1em}\mathop{\leqslant}\nolimits\hspace{0.1em}}
\def\ge{\hspace{0.1em}\mathop{\geqslant}\nolimits\hspace{0.1em}}

\mathchardef\inferieur="321E
\mathchardef\superieur="321F

\def\eqna{\begin{eqnarray*}}
\def\endeqna{\end{eqnarray*}}

\def\itemth#1{\item[${\mathrm{(#1)}}$]}

\def\verticalh{\vphantom{A^{\DS{A}}}}
\def\verticalhb{\vphantom{A^{\DS{A}}_{\DS{A}}}}
\def\verticalhbb{\DS{\vphantom{\frac{A^{\DS{A}}}{A_{\DS{A}}}}}}

\def\loewy{{\mathrm{LL}}}

\begin{document}

\title{Around Solomon's descent algebras}

\author{C. Bonnaf\'e \& G. Pfeiffer}
\address{\noindent 
Labo. de Math. de Besan\c{c}on (CNRS: UMR 6623), 
Universit\'e de Franche-Comt\'e, 16 Route de Gray, 25030 Besan\c{c}on
Cedex, France} 
\makeatletter
\email{bonnafe@math.univ-fcomte.fr}
\makeatother

\address{\noindent
Department of Mathematics,
National University of Ireland, Galway,
Ireland}
\email{goetz.pfeiffer@nuigalway.ie}

\subjclass{According to the 2000 classification:
Primary 20F55; Secondary 05E99}

\date{\today}

\begin{abstract} 
We study different problems related to the Solomon's descent algebra $\Sigma(W)$ 
of a finite Coxeter group $(W,S)$: positive elements, morphisms between 
descent algebras, Loewy length... One of the main result is that, if $W$ is 
irreducible and if the 
longest element is central, then the Loewy length of $\Sigma(W)$ is equal 
to $\DS{\Big\lceil\frac{|S|}{2}\Big\rceil}$.
\end{abstract}

\maketitle

\pagestyle{myheadings}

\markboth{\sc C. Bonnaf\'e \& G. Pfeiffer}{\sc Solomon's descent algebras}

\section*{Introduction}
Let $(W, S)$ be a finite Coxeter system.  The descent algebra
$\Sigma(W)$ of the finite Coxeter group $W$ is a subalgebra of the
group algebra $\QM W$ with a basis $\{x_I : I \subset S\}$, where $x_I$
is the sum in $\QM W$ of the distinguished coset representatives of the
parabolic subgroup $W_I$ in $W$.  It is a non-commutative preimage of
the ring of parabolic permutation characters of $W$, with respect to
the homomorphism $\theta$ which associates to $x_I$ the permutation
character of $W$ on the cosets of $W_I$.  Solomon~\cite{solomon}
discovered it as the real reason why the sign character of $W$ is a
linear combination of parabolic permutation characters.  He 
also showed that $\Ker \theta$ is the radical of $\Sigma(W)$.

The special case where $W$ is the symmetric group on $n$ points, i.e.,
a Coxeter group of type $A_{n-1}$, has received particular attention.
This type of descent algebra occurs as the dual of the Hopf algebra of
quasi-symmetric functions.  Atkinson~\cite{atkinson} has determined
the Loewy length of $\Sigma(W)$ in this case.

For general $W$, the descent algebra has been further studied as an
interesting object in its own right.  Bergeron, Bergeron, Howlett and
Taylor~\cite{bbht} have constructed explicit idempotents, decomposing
$\Sigma(W)$ into projective indecomposable modules.  Recently,
Blessenohl, Hohlweg and Schocker~\cite{BHS} could show that $\theta$
satisfies the remarkable symmetry $\theta(x)(y) = \theta(y)(x)$ for
all $x, y \in \Sigma(W)$.

The main purpose of this article is to determine the Loewy length of
$\Sigma(W)$ for all types of irreducible finite Coxeter groups $W$.
With the exception of type $D_n$, $n$ odd, this is done through a case
by case analysis, using computer calculations with 
{\sf CHEVIE}~\cite{chevie} for the exceptional types, in the final
Section~5.  Our results show in particular, that if $W$ is irreducible 
and if the longest
element $w_0$ is central in $W$ then the Loewy length of $\Sigma(W)$
is exactly $\DS{\Bigl\lceil\frac{|S|}2\Bigr\rceil}$, whereas in the other
cases, it lies between $\DS{\Bigl\lceil\frac{|S|}{2}\Bigr\rceil}$ and $|S|$.
Moreover, in Section~3, we study ideals generated by elements of
$\Sigma^+(W)$, the set of non-negative linear combinations of the
basis elements $x_I$ of $\Sigma(W)$, and show that the minimal
polynomial of an element of $\Sigma^+(W)$ is square-free.  Section~4
deals with various types of homomorphisms between descent algebras,
some restriction morphisms and one type related to self-opposed
subsets. A restriction morphism between the descent algebra of type 
$B_n$ and the descent algebra of type $D_n$ is also defined. 
Section~2 sets the scene in terms of a finite Coxeter group
$W$ and a length-preserving automorphism $\sigma$.  The general object
of interest is $\Sigma(W)^{\sigma}$, the subalgebra of fixed points of
$\sigma$ in $\Sigma(W)$.

\medskip

\noindent{\sc Remark - } 
If $W_n$ is a Weyl group of type $B_n$, there exists an extension 
$\Sigma_n'$ of the descent algebra $\Sigma(W_n)$ which 
was defined by Mantaci and Reutenauer \cite{mantaci} and 
studied by Hohlweg and the first author \cite{bh}. 
In \cite{bonnafe}, the first author investigates 
similar problems for this algebra (restriction morphisms, 
positive elements, Loewy series...): for instance, $\Sigma_n'$ 
has Loewy length $n$.

\medskip

\noindent{\sc Acknowledgement - } 
Some of the research leading to this paper was carried out when the
authors were visiting the Centre Interfacultaire Bernoulli at the
EPFL in Lausanne, Switzerland. They would like to express their 
gratitude for the Institute's hospitality. 

\bigskip

\section{Notation, preliminaries}

\medskip

\soussection{General notation} 
If $X$ is a set, $\PC(X)$ denotes the set of subsets of $X$ and 
$\PC^\#(X)$ denotes the set of proper subsets of $X$. If $k \in \ZM$, 
we denote by $\PC_{\le k}(X)$ the set of subsets $I$ of $X$ such that 
$|I| \le k$. 
The group algebra of a group $G$ over $\QM$ is denoted by $\QM G$. 
If $G$ is a finite group, let $\Irr G$ denote the set of its 
(ordinary) irreducible characters over $\CM$. 
The Grothendieck group of the category of finite dimensional 
$\CM G$-modules is identified naturally with the free $\ZM$-module 
$\ZM \Irr G$ and we set $\QM\Irr G=\QM \otimes_\ZM \ZM\Irr G$. 
If $A$ is a finite dimensional $\QM$-algebra, we denote by 
$\Rad A$ its radical. If $a \in A$, the centralizer of $a$ in $A$ 
is denoted by $Z_A(a)$. The set of irreducible characters of $A$ 
is denoted by $\Irr A$. 

\medskip

\soussection{Coxeter groups}
Let $(W,S)$ be a finite Coxeter group. Let $\ell : W \to \NM$ 
be the length function attached to $S$ and let $\le$ denote the 
Bruhat-Chevalley order on $W$. Let $w_0$ denote the longest 
element of $W$. If $I \in \PC(S)$, let $W_I$ denote the subgroup 
of $W$ generated by $I$. Recall that $(W_I,I)$ is a Coxeter group. 
The trivial character of $W_I$ is denoted by $1_I$. A {\it parabolic subgroup} 
of $W$ is a subgroup of $W$ which is conjugate to some $W_I$. 


\medskip

\soussection{Solomon descent algebra} 
If $I \subset S$, we set
$$
X_I=\{w \in W~|~\forall~s \in I,~ws > w\}. 
$$
Recall that an element $w \in W$ lies in $X_I$ 
if and only if $w(\D_I) \subset \Phi^+$. Let 
$$
x_I=\sum_{w \in X_I} w \in \QM W.
$$
Let 
$$\Sigma(W)=\mathop{\oplus}_{I \in \PC(S)} \QM x_I \quad \subset \quad\QM W.$$
If $\FC$ is a subset of $\PC(S)$, we set 
$$\Sigma_\FC(W)=\mathop{\oplus}_{I \in \FC} \QM x_I.$$
In particular, $\Sigma_{\PC(S)}(W)=\Sigma(W)$. 
Let $\th : \Sigma(W) \to \QM\Irr W$ be the unique linear map such that 
$\th(x_I)=\Ind_{W_I}^W 1_I$ 
for every $I \subset S$. Let $(\xi_I)_{I \in \PC(S)}$ denote the $\QM$-basis 
of $\Hom_\QM(\Sigma(W),\QM)$ dual to $(x_I)_{I \in \PC(S)}$. In 
other words,
$$x=\sum_{I \in \PC(S)} \xi_I(x) x_I$$
for every $x \in \Sigma(W)$. If $s \in S$, we write $x_s$ (resp. 
$\xi_s$) for $x_{\{s\}}$ (resp. $\xi_{\{s\}}$) for simplification.

If $I$ and $J$ are two subsets of $S$, we set 
$$
X_{I J}= (X_I)^{-1} \cap X_J.
$$
We write $I \equiv J$ if there exists $w \in W$ such that $J=\lexp{w}{I}$ 
(or, equivalently, if $W_I$ and $W_J$ are conjugate subgroups of $W$). 
The relation $\equiv$ is an equivalence relation on $\PC(S)$ and we denote by 
$\L$ the set of equivalence classes for this relation: it parametrizes 
the $W$-conjugacy classes of parabolic subgroups of $W$. 
We still denote by $\subset$ the order relation on $\L$ induced by inclusion. 
Let $\lamb : \PC(S) \to \L$ be the canonical surjection. 
We can now recall the following result of Solomon \cite{solomon}.

\bigskip

\noindent{\bf Solomon's Theorem.} 
{\it With the previous notation, we have:
\begin{itemize}
\itemth{a} If $I$ and $J$ are two subsets of $S$, then 
$$
x_I x_J = \sum_{d \in X_{I J}} x_{\lexp{d^{-1}}{I} \cap J}.
$$

\itemth{b} $\Sigma(W)$ is a unitary sub-$\QM$-algebra of $\QM W$. 

\itemth{c} $\th : \Sigma(W) \to \QM \Irr W$ is a morphism of 
$\QM$-algebras.

\itemth{d} $\Ker \th=\DS{\sum_{I \equiv J} \QM(x_I - x_J)}$. 

\itemth{e} $\Rad \Sigma(W)=\Ker \th$. 
\end{itemize}}

\bigskip

$\Sigma(W)$ is called {\it Solomon's descent algebra} of $W$. If $I$, $J$ and 
$K$ are three subsets of $S$, we set
$$X_{IJK}=\{d \in X_{IJ}~|~\lexp{d^{-1}}{I} \cap J = K\}.$$
Then, Solomon's Theorem (a) can be restated as follows:
\equat
x_I x_J = \sum_{K \in \PC(S)} |X_{IJK}|~x_K.
\endequat

\medskip

\soussection{Simple ${\boldsymbol{\Sigma(W)}}$-modules} 
The intersection of two parabolic subgroups of $W$ is a parabolic subgroup. 
Therefore, if $w \in W$, we define $W(w)$ to be the minimal 
parabolic subgroup of $W$ containing $w$. 
We denote by $\Lamb(w) \in \L$ the parameter of its conjugacy class. The map 
$$
\Lamb : W \longto \L
$$
is constant on conjugacy classes and is surjective: indeed, if $\l \in \L$, 
if $I \in \l$, and if $c$ is a Coxeter element of $W_I$, then $\Lamb(c)=\l$. 
The inverse image of $\l \in \L$ in $W$ is denoted by $\CC(\l)$. 
It is a union of conjugacy classes of $W$. 

If $\l \in \L$, let $\t_\l : \Sigma(W) \to \QM$, 
$x \mapsto  \th(x)(w)$, where $w \in \CC(\l)$. Recall that 
$\th(x)$ is a $\QM$-linear combination of permutation characters, so 
$\th(x)(w)$ lies in $\QM$. Moreover, $\t_\l$ does not depend on 
the choice of $w$ in $\CC(\l)$, and is a morphism of algebras. Also, the map 
$$
\fonction{\t}{\L}{\Irr \Sigma(W)}{\l}{\t_\l}
$$
is bijective. By definition, if $w \in W$ and $x \in \Sigma(W)$, then 
\equat\label{tau bete}
\t_{\Lamb(w)}(x)=\th(x)(w).
\endequat
Finally, recall that 
\equat\label{tau}
\t_{\lamb(J)}(x_I)=|X_{IJJ}|.
\endequat
It follows that 
\equat\label{tau bis}
xx_J - \t_{\lamb(J)}(x) x_J \in \Sigma_{\PC^\#(J)}(W)
\endequat
for every $x \in \Sigma(W)$. 

\section{Automorphisms of Coxeter groups}

\medskip

\soussection{General case}
We fix in this section an automorphism $\s$ of $W$ such that $\s(S)=S$. 
Since $\ell \circ \s=\ell$, $\s$ induces an automorphism of 
$\Sigma(W)$ which is still denoted by $\s$. 
The subalgebra of fixed points of $\s$ in $\Sigma(W)$ is denoted by $\Sigma(W)^\s$. 

\begin{lem}\label{radical sous-algebre}
Let $A$ be a sub-$\QM$-algebra of $\Sigma(W)$. Then $\Rad A = A \cap \Rad \Sigma(W)$. 
\end{lem}

\begin{proof}
Let $I=A \cap \Rad \Sigma(W)$. 
Since $\Sigma(W)$ is basic (i.e., all its simple modules 
are of dimension $1$), $\Rad \Sigma(W)$ is exactly the set 
of nilpotent elements of $\Sigma(W)$. Therefore, $\Rad A \subset \Rad \Sigma(W)$. 
In particular, $\Rad A \subset I$. Moreover, $I$ is a two-sided nilpotent 
ideal of $A$. So $I \subset \Rad A$ and we are done.
\end{proof}

\begin{coro}\label{radical fixe}
$\Rad\bigl(\Sigma(W)^\s\bigr)=\bigl(\Rad \Sigma(W)\bigr)^\s$.
\end{coro}

The automorphism $\s$ acts on $\PC(S)$ and this action induces an action 
of $\s$ on $\L$. The set of $\s$-orbits in $\L$ is denoted by $\L/\s$. 
It is easily checked that 
\equat\label{tau sigma}
\t_\l \circ \s^{-1} = \t_{\s(\l)}
\endequat
for every $\l \in \L$. In particular, if we denote by $\t_\l^\s$ the restriction 
of $\t_\l$ to $\Sigma(W)^\s$, then 
\equat\label{t}
\t_\l^\s=\t_{\s(\l)}^\s.
\endequat
It is also clear that $\t_\l^\s$ is an irreducible character of $\Sigma(W)^\s$. 

\begin{prop}\label{irreductibles fixe}
The map $\L \to \Irr\bigl(\Sigma(W)^\s\bigr)$, $\l \mapsto \t_\l^\s$ 
induces a bijection 
$\L/\s \simeq \Irr\bigl(\Sigma(W)^\s\bigr)$.
\end{prop}

\begin{proof}
By Corollary \ref{radical fixe} and since $\QM$ has characteristic $0$, 
$\th$ induces an isomorphism of algebras 
$$\Sigma(W)^\s/\Rad\bigl(\Sigma(W)^\s\bigr) \simeq (\im \th)^\s.$$
So we have a natural bijection between $\Irr \Sigma(W)^\s$ and 
$\Irr (\im \th)^\s$. 
If $\l \in \L$, let $e_\l$ be the idempotent of $\im \th$ 
such that $(\im \th)e_\l$ is a simple $\Sigma(W)$-module affording 
$\t_\l$. Then 
$$\im \th=\mathop{\oplus}_{\l \in \L} \QM e_\l$$ 
and $\s(e_\l)=e_{\s(\l)}$. So, 
$$(\im \th)^\s=\mathop{\oplus}_{\Omega \in \L/\s} \QM(\sum_{\l \in \O}e_\l).$$
This completes the proof of the proposition. 
\end{proof}

\medskip

\soussection{Action of ${\boldsymbol{w_0}}$} 
Let $\s_0$ denote the automorphism of $W$ induced by conjugation by $w_0$, 
the longest element of $W$. Then $\s_0(S)=S$, so $\s_0$ induces an 
automorphism of $\Sigma(W)$. Of course, we have
\equat
\Sigma(W)^{\s_0} = Z_{\Sigma(W)}(w_0).
\endequat
Let us introduce another classical basis of $\Sigma(W)$. If $w \in W$, we set
$$\RC(w)=\{s \in S~|~ws > w\}.$$
Then
\equat\label{r w0}
\RC(w_0w)=S \setminus \RC(w).
\endequat
If $J \in \PC(S)$, we set 
$$Y_J=\{w \in W~|~\RC(w)=J\}$$
$$y_J = \sum_{w \in Y_J} w \in \QM W.\leqno{\text{and}}$$
Then
\equat\label{xy}
x_I=\sum_{I \subset J} y_J,
\endequat
so $y_J \in \Sigma(W)$ and $(y_J)_{J \in \PC(S)}$ is a $\QM$-basis 
of $\Sigma(W)$. Note that $y_S=\{1\}$ and $y_\vide = \{w_0\}$, so $w_0 \in \Sigma(W)$. 
By \ref{r w0}, we have 
\equat\label{w0 y}
y_\vide y_J = w_0 y_J = y_{S \setminus J}.
\endequat

The centrality of $w_0$ can be characterized by the invertibilty of the
elements $y_J$.
\begin{prop}
  The longest element $w_0$ is central  in $W$ if and only if $y_J$ is
  invertible for all $J \in \PC(S)$.
\end{prop}

\begin{proof}
Clearly $x  \in \Sigma(W)$ is invertible  if and only  if $0 \not\in
\{\theta(x)(w) : w \in W\}$.  Moreover, $w_0 \not\in W_I$ unless $I = S$.
And by M\"obius inversion, 
$$y_J = \sum_{I \supset J} (-1)^{|I| - |J|} x_I
= x_S + \sum_{I \in \PC^\#(S): J \subset I} (-1)^{|I| - |J|} x_I.$$

Suppose $w_0$ is  central in $W$.  Then $w_0 \in  N_W(W_I)$ for all $I
\in  \PC(S)$ and  the index  $|N_W(W_I)  : W_I|$  is even  for $I  \in
\PC^\#(S)$.   Let  $w  \in  W$.   Then $\theta(x_I)(w)$,  which  is  a
multiple of  $|N_W(W_I) : W_I|$, is  even for $I  \in \PC^\#(S)$.  And
$\theta(y_J)(w)$, which is the  sum of $\pm \theta(x_I)(w)$ for certain
$I \in \PC^\#(S)$  and $\theta(x_S)(w) = 1$ is  odd, in particular not
zero.

Conversely, if $w_0$ is not central  in $W$, there is a maximal proper
subset $I \subset S$ such  that $I^{w_0} \neq I$.  (Otherwise $s^{w_0}
= s$ for all $s \in  S$, in contradiction to $w_0$ being non-central.)
It follows that, if  $w$ is an element of the same  shape as $I$, then
$\theta(x_I)(w) =  |N_W(W_I) :  W_I| =  1$.  Hence $y_I  = x_I  - x_S$
implies $\theta(y_I)(w) = 1 - 1 = 0$.
\end{proof}

If $I \in \PC(S)$, we set
$$x_I' = \sum_{K \in \PC(I)} \Bigl(-\frac{1}{2}\Bigr)^{|I|-|K|} x_K.$$
Note that $(x_I')_{I \in \PC(S)}$ is a basis of $\Sigma(W)$. Using \ref{xy}, 
it is easily checked that
\equat\label{xxy}
x_I' = \Bigl(-\frac{1}{2}\Bigr)^{|I|} \sum_{J \in \PC(S)} (-1)^{|I \cap J|} y_J.
\endequat
Therefore, by \ref{w0 y}, we get
\equat\label{w0 xx}
w_0 x_I' = (-1)^{|I|} x_I'.
\endequat
So, if $w_0$ is central in $W$, we can improve \ref{tau bis}:

\begin{lem}\label{w0 central}
Let $I \in \PC(S)$ and $x \in \Sigma(W)^{\s_0}$. Then
$$xx_I' \in \t_{\lamb(I)}(x) x_I' + \Sigma_{\PC_{\le |I|-2}(I)}(W).$$
\end{lem}

\begin{proof}
Let us write 
$$xx_I' = \sum_{J \subset I} \a_J x_J'.$$
Evaluating $\xi_I$ on each side, we get that $\a_I=\t_{\lamb(I)}(x)$ 
(see \ref{tau bis}). Since $x$ commutes with $w_0$, it follows 
from \ref{w0 xx} that $\a_J=0$ if $|J|-|I| \equiv 1 \mod 2$, 
as desired.
\end{proof}

\bigskip

\section{Positivity properties}

We denote by $\Sigma^+(W)$ the set of elements $a \in \Sigma(W)$ 
such that $\xi_I(a) \ge 0$ for every $I \in \PC(S)$. Note that 
$x_I \in \Sigma^+(W)$ for every $I \in \PC(S)$. If $a$, $b \in \Sigma^+(W)$, 
then 
\equat
a+b \in \Sigma^+(W)
\endequat
and, by Solomon's Theorem (a), 
\equat\label{positif produit}
ab \in \Sigma^+(W).
\endequat
The aim of this section is to study properties of the elements 
of $\Sigma^+(W)$ (ideals generated, minimal polynomial, 
centralizer...). 

\bigskip

\soussection{Ideals\label{section ideal}}
A subset $\FC$ of $\PC(S)$ is called {\it saturated} 
(resp. {\it equivariantly saturated})
if, for every $I \in \FC$ and every $I' \in \PC(S)$ such that 
$I' \subset I$ (resp. $\lamb(I')\subset \lamb(I)$), we have $I' \in \FC$. 
If $\FC$ is equivariantly saturated, then it is saturated. 
If $\FC$ is saturated (resp. equivariantly saturated) then, 
by Solomon's Theorem (a), $\Sigma_\FC(W)$ is a left (resp. two-sided) 
ideal of $\Sigma(W)$.

\medskip

\noindent{\sc Example and notation - } 
Then $\PC_{\le k}(S)$ is an equivariantly saturated subset of $\PC(S)$. 
Moreover, if $I \subset S$, then $\PC(I)$ and $\PC^\#(I)$ 
are saturated subsets of $\PC(S)$.\finl

\begin{prop}\label{caractere sature}
Let $\FC$ be a saturated subset of $\PC(S)$ and let $\chi_\FC$ 
denote the character of the left $\Sigma(W)$-module $\Sigma_\FC(W)$. Then 
$$
\chi_\FC = \sum_{I \in \FC} \t_{\lamb(I)}.
$$
\end{prop}

\begin{proof} 
This follows immediately from \ref{tau bis}.
\end{proof}

\medskip

If $a \in \Sigma(W)$, we set 
$$\FC(a)=\{I \in \PC(S)~|~\exists~J\in \PC(S),~
\bigl(\xi_J(a)\not=0\text{ and } I \subset J\bigr)\}$$
$$\FC_{\mathrm{eq}}(a)\{I \in \PC(S)~|~ \\
\exists~J\in \PC(S),~
\bigl(\xi_J(a)\not=0\text{ and } \lamb(I)\subset \lamb(J)\bigr)\}.$$
Note that $\FC(a) \subset \FC_{\mathrm{eq}}(a)$. 
Then $\FC(a)$ (resp. $\FC_{\mathrm{eq}}(a)$) is saturated (resp. 
equivariantly saturated) and, by Solomon's Theorem (a), 
\equat\label{gauche}
\Sigma(W)a\subset \Sigma_{\FC(a)}(W)
\endequat
and
\equat\label{droite}
a\Sigma(W)\subset \Sigma_{\FC_{\mathrm{eq}}(a)}(W).
\endequat
The next proposition shows that equality holds in \ref{droite} 
whenever $a \in \Sigma^+(W)$.

\begin{prop}\label{sature droite}
Let $a \in \Sigma^+(W)$. Then 
$$a\Sigma(W) = \Sigma_{\FC_{\mathrm{eq}}(a)}(W).$$
In particular, $\Sigma(W)a \subset a \Sigma(W)$. 
\end{prop}

\begin{proof}
We may, and we will, assume that $a\not= 0$. 
Let $\FC=\FC_{\mathrm{eq}}(a)$ and $\IC=a\Sigma(W)$. 
Then $\FC$ is equivariantly 
saturated and $\IC \subset \Sigma_\FC(W)$ (see \ref{droite}). 
Now let $I \in \FC$. We shall show by induction on $|I|$ that $x_I \in \IC$.

First, note that 
$$a x_\vide  = \Bigl(\sum_{I \in \PC(S)} |X_I| \xi_I(a) \Bigr) 
x_\vide$$
so, by hypothesis, $a x_\vide= m x_\vide$ with $m > 0$. Therefore, 
$x_\vide \in \IC$. Now, let $I \in \FC$ and assume that, for every 
$J \in \FC$ such that $|J| \le |I|-1$, we have $x_J \in \IC$. 
We want to prove that $x_I \in \IC$. 
Let $I_0 \in \PC(S)$ be such that $\lamb(I) \subset \lamb(I_0)$ and 
$\xi_{I_0}(a) \not= 0$. 
By the positivity of $a$ and by Solomon's Theorem (a), 
this shows that $\xi_I(a x_I) > 0$. 
But $a x_I = \sum_{J \in \PC(I)} \xi_J(a x_I) x_J$. 
Since $a x_I \in \IC$ and $x_J \in \IC$ for every $J \in \PC^\#(I)$, 
we get that $x_I \in \IC$, as desired.
\end{proof}

\begin{coro}\label{positif inversible}
Let $a \in \Sigma^+(W)$. Then $a$ is invertible 
in $\Sigma(W)$ if and only if $\xi_S(a) > 0$.
\end{coro}

%

\begin{coro}\label{somme ideaux}
Let $a_1$,\dots, $a_r \in \Sigma^+(W)$. Then $a_1 + \dots + a_r \in \Sigma^+(W)$ 
and 
$$a_1 \Sigma(W) + \dots + a_r \Sigma(W) = (a_1+\dots+a_r)\Sigma(W).$$
\end{coro}

\begin{proof}
By Proposition \ref{sature droite}, we have 
$$a_1 \Sigma(W) + \dots + a_r \Sigma(W) = 
\Sigma_{\FC_{\mathrm{eq}}(a_1) \cup \dots \cup \FC_{\mathrm{eq}}(a_r)}(W).$$
But it is clear that 
$$\FC_{\mathrm{eq}}(a_1) \cup \dots \cup \FC_{\mathrm{eq}}(a_r)=
\FC_{\mathrm{eq}}(a_1+\dots+a_r).$$
By applying Proposition \ref{sature droite} to $a=a_1+\dots+a_r$, 
we get the desired result.
\end{proof}

\bigskip

\soussection{Minimal polynomial} If $a \in \Sigma(W)$, we denote 
by $f_a(T) \in \QM[T]$ its minimal polynomial. Let $m_a : \Sigma(W) \to \Sigma(W)$, 
$x \mapsto ax$ be the left multiplication by $a$ and let $M_a$ be the matrix of 
$m_a$ in the basis $(x_J)_{J \in \PC(S)}$. 
The minimal polynomial of $a$ is equal to the minimal polynomial 
of the linear map $m_a$ (or of the matrix $M_a$). 
By \ref{tau bis}, $M_a$ is triangular 
(with respect to the order $\subset$ on $\PC(S)$) and its characteristic 
polynomial is 
$$\prod_{J \in \PC(S)} (T-\t_{\lamb(J)}(a)).$$
In particular
\equat
f_a\text{ \it is split over $\QM$.}
\endequat
The main result of this subsection is the following:

\begin{prop}\label{minimal positif}
Let $a \in \Sigma^+(W)$. Then $f_a$ is square-free.
\end{prop}

\begin{proof}
Before starting the proof, we gather in the next lemma some elementary properties 
of elements of $\Sigma^+(W)$. 
\begin{quotation}
\begin{lem}\label{croissance}
Let $I$, $J$ and $K$ be three subsets of $S$ such that $J \subset K$ and 
let $a \in \Sigma^+(W)$. Then:
\begin{itemize}
\itemth{a} $X_{IK} \subset X_{IJ}$ and $X_{IKK} \subset X_{IJJ}$. 

\itemth{b} $\t_{\lamb(K)}(a) \le \t_{\lamb(J)}(a)$.

\itemth{c} If $\t_{\lamb(K)}(a) = \t_{\lamb(J)}(a)$ and if 
$\xi_I(a) \not= 0$, then:
\begin{itemize}
\itemth{c1} $X_{IJJ}=X_{IKK}$.

\itemth{c2} If $J \varsubsetneq K$, then $X_{IKJ}=\vide$.
\end{itemize}
\end{itemize}
\end{lem}

\begin{proof}[Proof of Lemma \ref{croissance}]
It is clear that $X_{IK} \subset X_{IJ}$. Now, 
we have $X_{IKK}=\{d \in X_{IK}~|~K \subset \lexp{d^{-1}}{I}\}$, so (a) follows. 
Now, by \ref{tau}, we have 
$$\t_{\lamb(K)}(a)=\sum_{I \in \PC(S)} \xi_I(a) |X_{IKK}|.$$
So (b) and (c1) follow immediately from (a) and this equality. Let us now prove 
(c2). So assume that $\t_{\lamb(K)}(a) = \t_{\lamb(J)}(a)$, that 
$\xi_I(a) \not= 0$ and that $X_{IKJ}\neq\vide$. Let $d \in X_{IKJ}$. 
Then $d \in X_{IJ}$ by (a) and 
$J =\lexp{d^{-1}}{I} \cap K \subset \lexp{d^{-1}}{I}$. 
In other words, $d \in X_{IJJ}$. Therefore, $d \in X_{IKK}$ by (c1) and, since 
$J =\lexp{d^{-1}}{I} \cap K$, we have $J=K$, as expected.
\end{proof}
\end{quotation}

Let $\xi \in \QM$ be an eigenvalue of $m_a$. Let 
$\FC=\{J \in \PC(S)~|~\t_{\lamb(J)}(a)=\xi\}$. Note that $\FC\not= \vide$. 
Since the matrix $M_a=(\xi_J(ax_K))_{K,J \in \PC(S)}$ is triangular, 
it is sufficient to show that 
the square matrix $(\xi_J(ax_K))_{K,J \in \FC}$ is diagonal. So, let 
$J$ and $K$ be two elements of $\FC$ such that $\xi_J(ax_K) \not= 0$. 
We want to show that $J=K$. 
First, since $\xi_J(ax_K) \not= 0$, we have $J \subset K$. Moreover, 
there exists $I \in \PC(S)$ such that $\xi_I(a) \not= 0$ and $X_{IKJ} \not= \vide$. 
But, by (2), we have $X_{IJJ}=X_{IKK}$. Now, let $d \in X_{IKJ}$ (such 
a $d$ exists by hypothesis). Then $d \in X_{IJ}$ and 
$J =\lexp{d^{-1}}{I} \cap K \subset \lexp{d^{-1}}{I}$. 
In other words, $d \in X_{IJJ}$. Therefore, $d \in X_{IKK}$ and, since 
$J =\lexp{d^{-1}}{I} \cap K$, we have $J=K$, as expected.
\end{proof}

\begin{coro}\label{puissance ideal}
Let $a \in \Sigma^+(W)$ and let $n \ge 1$. Then $a^n\Sigma(W)=a \Sigma(W)$ 
and $\Sigma(W)a^n=\Sigma(W)a$.
\end{coro}

\begin{proof}
It is sufficient to prove this result for $n=2$. If $a$ is invertible, then 
the result is clear. If $a$ is not invertible then, by Proposition 
\ref{minimal positif}, the minimal polynomial $f_a$ of $a$ is divisible 
by $T$ and not by $T^2$. This shows that $a \in \QM[a] a^2 = \QM[a] a^2$. 
So $a^2 \in \Sigma(W) a$ and $a^2 \in a \Sigma(W)$, as expected.
\end{proof}

\begin{coro}\label{espace propre}
Let $M$ be a $\Sigma(W)$-module and let $\chi_M$ denote its character. 
Write $\chi_M=\t_{\l_1} + \dots + \t_{\l_r}$, with $\l_1$,\dots, $\l_r \in \L$ 
(possibly non-distinct). Let $a \in \Sigma^+(W)$ and let $\xi \in \QM$. 
Then
$$\dim_\QM \Ker(a-\xi \Id_M~|~M) = |\{1 \le i \le r~|~\t_{\l_i}(a)=\xi\}|.$$
\end{coro}

\begin{proof}
Indeed, if $x \in \Sigma(W)$, then $(\t_{\l_i}(x))_{1 \le i \le r}$ 
is the multiset of eigenvalues of $x$ in its action on $M$. 
But, by Proposition \ref{minimal positif}, $a$ acts semisimply on $M$. 
This proves the result.
\end{proof}

\example{reguliere} 
Consider here the left $\Sigma(W)$-module $\QM W$, with the natural 
action by left multiplication. Let $\chi$ denote its character. Then it is 
easy and well-known that 
$$\chi(x_I)=|W|$$
for every $I \in \PC(S)$. Therefore, 
$$\chi=\sum_{\l \in \L} |\CC(\l)|~\t_\l = \sum_{w \in W} \t_{\Lamb(w)}.\leqno{(\arm)}$$
Indeed, by \ref{tau bete}, we have 
$$\sum_{w \in W} \t_{\Lamb(w)}(x_I) = \sum_{w \in W} \th(x_I)(w) = 
|W| \langle \th(x_I), 1_S\rangle = |W|.$$
Therefore, if $a \in \Sigma^+(W)$ and $\xi \in \QM$, it follows from 
Corollary \ref{espace propre} and \ref{tau bete} that 
$$\dim_\QM \Ker(a-\xi \Id_{\QM W}~|~\QM W) = |\{w \in W~|~\th(a)(w)=\xi\}|.~
\SS{\square}\leqno{(\brm)}$$

\bigskip

\soussection{Centralizers}
The aim of this subsection is to prove 
a few results on the dimension of the centralizer of elements of 
$\Sigma^+(W)$. We first start with some easy observation. 

Let $a \in \Sigma(W)$. 
Let $\mu_a : \Sigma(W) \to \Sigma(W)$, $x \mapsto ax -xa$. Then 
\equat\label{noyau centralise}
\Ker \mu_a = Z_{\Sigma(W)}(a)
\endequat 
so that
\equat\label{dif dim}
\dim_\QM \Sigma(W) - \dim_\QM Z_{\Sigma(W)}(a) = \dim_\QM(\im \mu_a).
\endequat

\begin{prop}\label{dimension bornee}
Let $a \in \Sigma^+(W)$. 
Then
$$\dim_\QM Z_{\Sigma(W)}(a) = 2^{|S|} -
|\FC_{\mathrm{eq}}(a)| + \dim_\QM \Sigma(W) a - \dim_\QM(\im \mu_a \cap \Sigma(W)a).$$
In particular, 
$$\dim_\QM Z_{\Sigma(W)}(a) \le 
2^{|S|}-|\FC_{\mathrm{eq}}(a)| + \dim_\QM \Sigma(W) a
\le  2^{|S|}-|\FC_{\mathrm{eq}}(a)| + |\FC(a)|.$$
\end{prop}

\noindent{\sc Remark - } Recall that $\FC(a) \subset \FC_{\mathrm{eq}}(a)$ 
so that the right-hand side 
of the above inequality is always $\le  2^{|S|}$.\finl

\begin{proof}
Let $\pi : a\Sigma(W) \to a\Sigma(W)/\Sigma(W)a$ be the 
canonical projection. Let $f : \Sigma(W) \to a\Sigma(W)$, $x \mapsto ax$. 
Then $\pi \circ f$ is surjective by definition. Moreover, note that 
the image of $\mu_a$ is contained in $a\Sigma(W)$. By 
definition, $\pi \circ f=\pi \circ \mu_a$. Therefore, $\pi \circ \mu_a$ is surjective. 
In particular, 
$$\dim_\QM(\im \mu_a) = |\FC_{\mathrm{eq}}(a)| - \dim_\QM \Sigma(W) a 
+ \dim_\QM(\im \mu_a \cap \Sigma(W)a).$$
So the result now follows from \ref{dif dim}.
\end{proof}

\example{strict}
The following example shows that the first inequality in Proposition 
\ref{dimension bornee} might be strict. Assume here that 
$W=\SG_4$ is of type $A_3$. Write $S=\{s_1,s_2,s_3\}$, 
with $s_1s_3=s_3s_1$. Then 
$$\dim_\QM Z_{\Sigma(W)}(x_{\{s_1,s_2\}}) = 5$$
$$2^{|S|}-|\FC_{\mathrm{eq}}(x_{\{s_1,s_2\}})|+\dim_\QM \Sigma(W)x_{\{s_1,s_2\}} 
= 6.~\SS{\square}\leqno{\text{and}}$$

\begin{coro}\label{dimension bornee coro}
Let $a \in \Sigma^+(W)$ and assume that $\Rad \Sigma(W) \cap \Sigma(W) a = 0$. Then 
$$\dim_\QM Z_{\Sigma(W)}(a) = 
2^{|S|}-|\FC_{\mathrm{eq}}(a)| + \dim_\QM \Sigma(W) a.$$
\end{coro}

\begin{proof}
Since $\im \mu_a \subset \Rad \Sigma(W)$, the hypothesis implies that 
$\im \mu_a \cap \Sigma(W) a =0$. So the result follows now from 
Proposition \ref{dimension bornee}.
\end{proof}

\example{xs}
Let $s \in S$. 
Let $C(s)$ denote the set of elements of $S$ which are conjugate to 
$s$ in $W$ and let $c(s)=|C(s)|$. Let 
$a=\sum_{t \in C(s)} \a_t x_t \not= 0$ be such that $\a_t \ge 0$ for every 
$t \in C(s)$. Then 
$$\FC_{\mathrm{eq}}(a)=C(s) \cup \{\vide\}.\leqno{(1)}$$
Moreover, if $t \in C(s)$, then $\xi_s(xx_s)=\xi_t(xx_t)$ 
for every $x \in \Sigma(W)$ (see \ref{tau} and \ref{tau bis}). Therefore, 
$$\Sigma(W) a = \QM a \oplus \QM x_\vide.\leqno{(2)}$$
In particular, 
$$\Rad \Sigma(W) \cap \Sigma(W) a =0.\leqno{(3)}$$
It then follows from (1), (2), (3) and Corollary \ref{dimension bornee coro} that 
$$\dim_\QM Z_{\Sigma(W)}(a)= 2^{|S|} - c(s)+1.\leqno{(4)}$$
Note that this equality holds if $a=x_s$.\finl

\begin{coro}\label{puissance}
Let $a \in \Sigma^+(W)$ and let $n \ge 1$. Then
$Z_{\Sigma(W)}(a^n)=Z_{\Sigma(W)}(a)$. 
\end{coro}

\begin{proof}
Since $Z_{\Sigma(W)}(a^n) \subset Z_{\Sigma(W)}(a)$, we only need to prove 
that the dimensions of both centralizers are equal. 
First, by Corollary \ref{puissance ideal}, we have 
$\dim_\QM \Sigma(W)a = \dim_\QM \Sigma(W)a^n$ and 
$\dim_\QM a\Sigma(W) = \dim_\QM a^n\Sigma(W)$. 
So, by Proposition \ref{dimension bornee}, we only need to prove that 
$$\dim_\QM(\im \mu_{a^n} \cap \Sigma(W) a) = \dim_\QM(\im \mu_a \cap \Sigma(W)a).
\leqno{(P_n)}$$
We will show $(P_n)$ by induction on $n$, the case where $n=1$ being 
trivial. So we assume that $n \ge 2$ and that 
$(P_{n-1})$ holds. 
First, note that $\mu_{a^n}(x)=a \mu_{a^{n-1}}(x) + \mu_a(x) a^{n-1}$. 
Therefore, $\mu_{a^n}(x) \in \Sigma(W)a$ if and only if 
$a \mu_{a^{n-1}}(x) \in \Sigma(W)a$. But, by Corollary 
\ref{puissance ideal}, the map $\kappa : a \Sigma(W) \to a \Sigma(W)$, 
$u \mapsto au$ is an isomorphism and stabilizes $\Sigma(W)a$. Therefore, 
$a \mu_{a^{n-1}}(x) \in \Sigma(W)a$ if and only if $\mu_{a^{n-1}}(x) \in \Sigma(W)a$. 
In other words, 
$$\im \mu_{a^n} \cap \Sigma(W) a = \kappa^{-1}(\im \mu_{a^{n-1}} \cap \Sigma(W) a).$$
This shows $(P_n)$.
\end{proof}

\bigskip

\soussection{Counter-examples} 
In this subsection, we provide examples to show that the different 
results of this section might fail if the positivity property 
is not satisfied. 

\medskip

$\bullet$ {\sc First statement of Proposition \ref{sature droite} - } 
Assume here that $W=\SG_3$ is of type $A_2$ and write $S=\{s_1,s_2\}$. 
Let $a=x_{s_1}-x_{s_2}$. Then $\Rad \Sigma(W) = \QM a$. 
Therefore, $a \Sigma(W) = \QM a \not= \Sigma_{\FC_{\mathrm{eq}}(a)}(W)$. 

\medskip

$\bullet$ {\sc Second statement of Proposition \ref{sature droite} - } 
Assume here that $W=\SG_4$ is of type $A_3$. 
Write $S=\{s_1,s_2,s_3\}$, with $s_1s_3=s_3s_1$. 
Let $a=x_{s_1}-x_{s_2,s_3}$. Then 
$x_{s_2}-x_{s_3}$ belongs to $\Sigma(W)a$ but does 
not belong to $a \Sigma(W)$.

\medskip

$\bullet$ {\sc Corollary \ref{positif inversible} - } 
Assume here that $W=\SG_3$ is of type $A_2$ and write $S=\{s_1,s_2\}$. 
Let $a=x_S -x_{s_2}$. Then $\xi_S(a) > 0$ but $a$ is not invertible. 

\medskip

$\bullet$ {\sc Corollary \ref{somme ideaux} - } 
Let $a \in \Sigma^+(W)$ be non-zero. Then 
$a\Sigma(W) + (-a)\Sigma(W)=a\Sigma(W) \not= (a + (-a))\Sigma(W)=0$.

\medskip

$\bullet$ {\sc Proposition \ref{minimal positif} and Corollary \ref{puissance ideal} - } 
Let $a \in \Rad \Sigma(W)$ be non-zero. Then $f_a(T)=T^n$ 
for some $n \ge 2$, so $f_a(T)$ is not square-free. 
Moreover, $\Sigma(W)a \not= \Sigma(W) a^n=0$ and 
$a\Sigma(W) \not= a^n\Sigma(W)=0$.

\medskip

$\bullet$ {\sc Corollary \ref{puissance} - } 
Let $a \in \Rad \Sigma(W)$ be non-central. Then there exists $n \ge 2$ such that 
$a^n=0$ is central.

\bigskip

\section{Some morphisms between Solomon descent algebras}

\medskip

\soussection{Restriction morphisms} 
Whenever $K \subset S$, F. Bergeron, N. Bergeron, R.B. Howlett 
and D.E. Taylor have constructed a so-called {\it restriction morphism} 
between $\Sigma(W)$ and $\Sigma(W_K)$. They do not say that they are 
morphisms of algebras but this can be deduced from 
some of their results \cite[13 and Proposition 2.6]{bbht}. However, 
we present here a simpler proof (see Proposition \ref{bbht}). 

In this subsection, we recall 
the definition and the basic properties of these restriction 
morphisms, and we prove some results on their image. We first 
need some notation:

\medskip

\begin{quotation}
\noindent{\sc Notation - } If $K \subset S$, we denote by $X_I^K$, $x_I^K$, $\th_K$, 
$\equiv_K$, $\L_K$, $\lamb_K$ and $\t_\l^K$ for the objects defined in $W_K$ 
instead of $W$ and which correspond respectively to 
$X_I$, $x_I$, $\th$, $\equiv$, $\L$, $\lamb$ and $\t_\l$.
\end{quotation}

\medskip

If $K \subset S$, let 
$\Res_K : \Sigma(W) \to \Sigma(W_K)$ 
denote the $\QM$-linear map such that 
$$
\Res_K(x_I)= \sum_{d \in X_{KI}} x_{K \cap \lexp{d}{I}}^K
$$
for every $I \in \PC(S)$. If $K \subset L \subset S$, we denote 
by $\Res_K^L : \Sigma(W_L) \to \Sigma(W_K)$ the map 
defined like $\Res_K$ but inside $W_L$. Finally, if $K' \in \PC(S)$ 
and if $d \in X_{KK'}$ are such that $\lexp{d}{K'}=K$, then 
the map $d_* : \Sigma(W_{K'}) \to \Sigma(W_K)$, $x \mapsto dxd^{-1}$ 
is well-defined and is an isomorphism of algebras. It sends 
$x_I^{K'}$ to $x_{\lexp{d}{I}}^K$ ($I \in \PC(K')$). Let us gather 
in the next proposition the formal properties of the map $\Res_K$:

\begin{prop}\label{bbht}
Let $K \in \PC(S)$. Then:
\begin{itemize}
\itemth{a} If $x \in \Sigma(W)$, then $x_K \Res_K(x)=xx_K$.

\itemth{b} $\Res_K$ is an homomorphism of algebras.

\itemth{c} If $K \subset L \subset S$, then 
$\Res_K = \Res_K^L \circ \Res_L$.

\itemth{d} The diagram
$$\diagram
\Sigma(W) \rrto^{\DS{\th}} \ddto_{\DS{\Res_K}} && 
\QM\Irr W \ddto^{\DS{\Res_{W_K}^W}} \\
&&\\
\Sigma(W_K) \rrto^{\DS{\th_K}} && \QM\Irr W_K \\
\enddiagram$$
is commutative.

\itemth{e} If $K' \in \PC(S)$ 
and if $d \in X_{KK'}$ are such that $\lexp{d}{K'}=K$, then 
$$\Res_K = d_* \circ \Res_{K'}.$$
\end{itemize}
\end{prop}

\begin{proof}
If $I \subset K$, then $x_K x_I^K=x_I$. So the map 
$\mu_K : \Sigma(W_K) \to \Sigma(W)$, $x \mapsto x_K x$ is well-defined 
and injective. 
(a) follows from this observation and from Solomon's Theorem (a).
(b) and (c) follow from (a) and from the injectivity of $\mu_K$ 
(note that $\Res_K(1)=1$). 
(d) is a direct consequence of the Mackey formula.
Finally, we have $x_{K'}=x_K d$. So (e) follows again from (a) 
and from the injectivity of $\mu_K$. 
\end{proof}

\bigskip

The natural map $\PC(K) \to \PC(S)$ induces a map 
$\pi_K : \L_K \to \L$. The next corollary is a generalisation of 
\cite[Theorem 3.6]{atkinson}.

\begin{coro}\label{factorisation simples}
If $\l \in \L_K$, then $\t_{\pi_K(\l)} = \t_\l^K \circ \Res_K$.
\end{coro}

\begin{proof}
This follows from Proposition \ref{bbht} (d).
\end{proof}

The Corollary \ref{factorisation simples} can be written as follows: if 
$I \in \PC(K)$, then
\equat\label{factorisation bis}
\t_{\lamb(I)} = \t_{\lamb_K(I)} \circ \Res_K.
\endequat
The next result generalizes \cite[Theorem 2.3]{atkinson}.

\begin{prop}\label{decomposition}
$\Sigma(W) = \Ker \Res_K \oplus \Sigma(W) x_K$ and 
$\Ker \Res_K$ is the set of $x \in \Sigma(W)$ such that 
$xx_K=0$. 
\end{prop}

\begin{proof}
By Proposition \ref{bbht} (a), we have 
$\dim_\QM(\im \Res_K) = \dim_\QM \Sigma(W)x_K$. 
Therefore,
$$\dim_\QM(\Ker \Res_K) + \dim_\QM \Sigma(W) x_K = 
\dim_\QM \Sigma(W).$$
Now, let $x \in \Sigma(W)$ be such that $xx_K \in \Ker \Res_K$. 
According to the previous equality, it is sufficient to show that $xx_K = 0$. 
But, by Proposition \ref{bbht} (a), we have that $xx_K^2=0$. 
By Corollary \ref{puissance ideal}, this implies that $xx_K=0$.
\end{proof}

\begin{coro}\label{test}
The following are equivalent:
\begin{itemize}
\itemth{1} $\Res_K$ is surjective.

\itemth{2} $\dim \Sigma(W) x_K = 2^{|K|}$.

\itemth{3} $\Sigma(W) x_K = \Sigma_{\PC(K)}(W)$.
\end{itemize}
\end{coro}

\begin{proof}
By Proposition \ref{bbht} (a), we have that 
$\dim_\QM \Sigma(W)x_K = \dim_\QM(\im \Res_K)$. By Solomon's Theorem (a), 
we have that $\Sigma(W)x_K \subset \Sigma_{\PC(K)}(W)$. Moreover, note that 
$\dim_\QM \Sigma_{\PC(K)}(W)=2^{|K|}$. The corollary follows from 
Proposition \ref{decomposition} and these 
three observations. 
\end{proof}

\bigskip

We now investigate further the image of $\Res_K$. 
First, let
$$W(K)=\{w \in X_{KK}~|~\lexp{w}{K}=K\}.$$
Then $W(K)$ is a subgroup of $W$ and 
$$N_W(W_K)= W(K) \ltimes W_K.$$
Moreover, $W(K)$ acts on $\Sigma(W_K)$ by conjugation. 

\begin{prop}\label{points fixes}
$\im \Res_K \subset \Sigma(W_K)^{W(K)}$.
\end{prop}

\begin{proof}
This follows immediately from Proposition \ref{bbht} (e).
\end{proof}

\begin{coro}\label{atkinson surjectif}
If $\Res_K$ is surjective, then 
the map $\pi_K : \L_K \to \L$ is injective and $W(K)$ acts trivially 
on $W_K$. 
\end{coro}

\begin{proof}
This follows immediately from Corollary \ref{factorisation simples} 
and Proposition \ref{points fixes}.
\end{proof}

\example{type A surjectif} 
Assume here that $W=\SG_n$ is the symmetric group of degree $n$. 
View $\SG_{n-1}$ as a parabolic subgroup of $W$. Then, 
by \cite{BGR}, the restriction morphism $\Sigma(\SG_n) \to \Sigma(\SG_{n-1})$ 
is surjective. Therefore, by Proposition \ref{bbht} (c) and (e), 
if $K$ is a subset of $S$ such that $W_K$ is irreducible, 
then $\Res_K$ is surjective. 

Moreover, the map $\pi_K$ is injective if and only if $W_K$ is irreducible. 
So we have shown that, if $W$ is irreducible of type $A$, then 
$\Res_K$ is surjective if and only if $\pi_K$ is injective.\finl

\medskip

\examples{surjectif exceptionnels}
Let $W$ be irreducible of exceptional type. Let $n=|S|$. We 
write $S=\{s_1,s_2,\dots,s_n\}$ following the convention of Bourbaki 
\cite[Planches I-IX]{bourbaki}. For simplification, we denote 
by $i_1i_2\dots i_k$ the subset $\{s_{i_1},s_{i_2},\dots,s_{i_k}\}$ of $S$ 
(for instance, $134$ stands for $\{s_1,s_3,s_4\}$). Then, computations 
using {\tt CHEVIE} show that:
\begin{itemize}
\itemth{a} If $W$ is of type $E_6$, $E_7$, $E_8$, $G_2$ or $H_3$, then 
$\Res_K$ is surjective if and only if $|K|\in \{0,1,|S|\}$.

\itemth{b} If $W$ is of type $F_4$, then $\Res_K$ is surjective 
if and only if $K$ belongs to $\{1234,123,234,13,14,23,24,1,2,3,4,\vide\}$. 

\itemth{c} If $W$ is of type $H_4$, then $\Res_K$ is surjective 
if and only if $K$ belongs to $\{1234,123,1,2,3,4,\vide\}$.\finl 
\end{itemize}

\bigskip

\noindent{\sc Remark - } 
The examples \ref{surjectif exceptionnels} show that the converse 
of Corollary \ref{atkinson surjectif} is not true in general.\finl

\bigskip

We will see in the next subsection some other examples of restriction 
morphisms (groups of type $B$ or $D$) and some results concerning their images. 

\bigskip

\soussection{Type ${\boldsymbol{B}}$, type ${\boldsymbol{D}}$: another restriction 
morphism\label{section BD}}
We shall investigate here some properties of $\Sigma(W)$ whenever 
$W$ is of type $B$ or $D$. We fix in this subsection a natural number $n \ge 1$. 
Let $(W_n,S_n)$ be a Coxeter group of type $B_n$. We write 
$S_n=\{t,s_1,s_2,\dots,s_{n-1}\}$ in such a way that the Dynkin 
diagram of $W_n$ is 
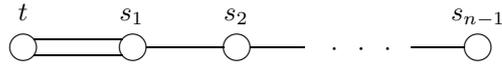
\begin{center}
\begin{picture}(220,30)
\put( 40, 10){\circle{10}}
\put( 44,  7){\line(1,0){33}}
\put( 44, 13){\line(1,0){33}}
\put( 81, 10){\circle{10}}
\put( 86, 10){\line(1,0){29}}
\put(120, 10){\circle{10}}
\put(125, 10){\line(1,0){20}}
\put(155,  7){$\cdot$}
\put(165,  7){$\cdot$}
\put(175,  7){$\cdot$}
\put(185, 10){\line(1,0){20}}
\put(210, 10){\circle{10}}
\put( 38, 20){$t$}
\put( 76, 20){$s_1$}
\put(115, 20){$s_2$}
\put(200, 20){$s_{n{-}1}$}
\end{picture}
\end{center}
Let $s_1'=ts_1t$, $S_n'=\{s_1',s_1,s_2,\dots,s_{n-1}\}$ and $W_n'=<S_n'>$. 
Then $(W_n',S_n')$ is a Weyl group of type $D_n$: its Dynkin diagram is 
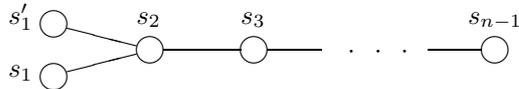
\begin{center}
\begin{picture}(220,40)
\put( 45, 25){\circle{10}}
\put( 45,  5){\circle{10}}
\put( 50,  23){\line(4,-1){26}}
\put( 50, 7){\line(4,1){26}}
\put( 81, 15){\circle{10}}
\put( 86, 15){\line(1,0){29}}
\put(120, 15){\circle{10}}
\put(125, 15){\line(1,0){20}}
\put(155, 12){$\cdot$}
\put(165, 12){$\cdot$}
\put(175, 12){$\cdot$}
\put(185, 15){\line(1,0){20}}
\put(210, 15){\circle{10}}
\put(28,5){$s_1$}
\put(28,25){$s_1'$}
\put( 76, 25){$s_2$}
\put(115, 25){$s_3$}
\put(200, 25){$s_{n{-}1}$}
\end{picture}
\end{center}
Recall that $W_n = <t> \ltimes W_n'$. 
So $X_n=\{1,t\}$ is the set of minimal length coset representatives 
of $W_n/W_n'$. We set $x_n=1+t \in \QM W_n$. 
Note that conjugacy by $t$ induces the unique 
non-trivial automorphism of $W_n'$ which stabilizes $S_n'$: this automorphism 
will be denoted by $\s_n$. If $I \subset S_n'$ or if $I \subset S_n$, we denote 
by $W_I$ the subgroup of $W_n$ generated by $I$. It is a standard parabolic 
subgroup of $W_n'$ or of $W_n$ and it might be a parabolic subgroup of both. 
If $I \subset S_n'$, we still denote by $X_I^{S_n}$ the set of $w \in W_n$ 
such that $w$ has minimal length in $w W_I$ and we set 
$x_I^{S_n}=\sum_{w \in X_I^{S_n}} w \in \QM W_n$. Therefore, 
if $I \subset S_n'$, 
\equat\label{?}
x_I^{S_n}=(1+t) x_I^{S_n'}.
\endequat

If $I \subset S_n$, then it is easy to check that 
\equat\label{conjugaison BD}
W_I \cap W_n' = W_{W_I \cap S_n'}\qquad\text{and}\qquad
\lexp{t}{W_I} \cap W_n' = W_{\lexp{t}{W_I} \cap S_n'}.
\endequat
Moreover, if $t \not\in I$, then 
\equat\label{conjugaison BD bis}
W_I \cap W_n' = W_I\qquad\text{and}\qquad
\lexp{t}{W_I} \cap W_n' = W_{\lexp{t}{I}}.
\endequat
We set
$$X_{I,n}=X_n \cap X_I^{-1}$$
and
$$\Res_n x_I^{S_n} = \sum_{d \in X_{I,n}} x_{\lexp{d^{-1}}{W_I} \cap S_n'}^{S_n'} \in 
\Sigma(W_n').$$
In other words, by \ref{conjugaison BD} and \ref{conjugaison BD bis}, 
\equat\label{def res}
\Res_n x_I^{S_n}=\begin{cases}
x_{W_I \cap S_n'}^{S_n'} & \text{if $t \in I$,} \\
x_I^{S_n'} + x_{\lexp{t}{I}}^{S_n'}& \text{if $t \not\in I$.} \\
\end{cases}
\endequat

This can be extended by linearity to a map $\Res_n : \Sigma(W_n) \to \Sigma(W_n')$. 
This map shares with the restriction morphisms many properties:

\begin{prop}\label{restriction BD}
With the above notation, we have:
\begin{itemize}
\itemth{a} If $x \in \Sigma(W_n)$, then $x_n \Res_n(x)=xx_n$.

\itemth{b} $\Res_n$ is an homomorphism of algebras.

\itemth{c} $\Res_{S_{n-1}'}^{S_n'} \circ \Res_n = \Res_{n-1} \circ \Res_{S_{n-1}}^{S_n}$.

\itemth{d} The diagram
$$\diagram
\Sigma(W_n) \rrto^{\DS{\th_n}} \ddto_{\DS{\Res_n}} && 
\QM\Irr W_n \ddto^{\DS{\Res_{W_n'}^{W_n}}} \\
&&\\
\Sigma(W_n') \rrto^{\DS{\th_n'}} && \QM\Irr W_n' \\
\enddiagram$$
is commutative.

\itemth{e} $\im \Res_n = \Sigma(W_n')^{\s_n}$.
\end{itemize}
\end{prop}

\begin{proof}
(a) Let $I \subset S_n$. We want to prove that $x_I^{S_n} (1+t)=(1+t) \Res_n(x_I^{S_n'})$. 
First, assume that $t \not\in I$. Then $W_I \subset W_n'$. Therefore, 
$x_I^{S_n}=(1+t) x_I^{S_n'}$. Consequently,
\eqna
x_I^{S_n} (1+t)&=&(1+t) x_I^{S_n'} (1+t) \\
&=& x_I^{S_n'}+t x_I^{S_n'} + x_I^{S_n'} t + tx_I^{S_n'} t \\
&=& (1+t)(x_I^{S_n'} + x_{\lexp{t}{I}}^{S_n'}) \\
&=& (1+t) \Res_n(x_I^{S_n}),
\endeqna
as desired. Now, assume that $t \in I$. Then $X_n=\{1,t\}$ is a set 
of minimal length coset representatives of $W_I/(W_I \cap W_n')$. 
Therefore $x_I^{S_n} x_n = x_{W_I \cap S_n'}^{S_n}=x_n x_{W_I \cap S_n'}^{S_n'}$, 
as expected (note that the last equality follows from \ref{?}. 

\medskip

(b) First, note that $\Res_n(1)=\Res_n(x_{S_n}^{S_n})=x_{S_n'}^{S_n'}=1$ by 
definition. The fact that $\Res_n(xy)=\Res_n(x)\Res_n(y)$ 
for all $x$, $y \in \Sigma(W_n)$ follows immediately from (a) and 
from the fact that the map $\mu_n : \QM W_n' \to \QM W_n$, $x \mapsto x_n x$ is 
injective. 

\medskip

(c) follows also from (a) and from the fact 
$x_n x_{S_{n-1}'}^{S_n'}=x_{S_{n-1}}^{S_n} x_{n-1}$. 

\medskip

(d) follows from the Mackey formula for tensor product of induced 
characters. 

\medskip

(e) This follows easily from \ref{def res}.
\end{proof}

\bigskip

We conclude this subsection by two examples where the image of the 
restriction map $\Res_K$ is computed explicitly. The first one concerns 
type $B$ (see Proposition \ref{surjectif B}) while the second 
one concerns the type $D$ (see Corollary \ref{surjectif D}). 

\begin{prop}\label{surjectif B}
The map $\Res_{S_{n-1}}^{S_n} : \Sigma(W_n) \to \Sigma(W_{n-1})$ 
is surjective.
\end{prop}

\begin{proof}
We have
\begin{multline*}
X_{S_{n-1}}^{S_n} = \{ s_i s_{i+1} \dots s_{n-1}~|~1 \le i \le n\} \\
\coprod \quad\{s_i s_{i-1} \dots s_1 ts_1s_2 \dots s_{n-1}~|~0 \le i \le n-1\}.
\end{multline*}
Therefore, if $d \in W_n$ and if $i \in \{1,2,\dots,n-1\}$ are such that 
$d^{-1} \in X_{S_{n-1}}^{S_n}$, $ds_i > d$, and $d s_i d^{-1} \in S_{n-1}$, then 
$$ds_i d^{-1} \in \{s_i,s_{i-1}\}.\leqno{(*)}$$
We define a total order $\infspe$ on $\PC(S_{n-1})$. Let $I$ and $J$ 
be two subsets of $S_{n-1}$. Then we write $I \infspe J$ 
if and only if one of the following two conditions are satisfied:
\begin{quotation}
\begin{itemize}
\itemth{1} $|I| < |J|$

\itemth{2} $|I|=|J|$ and $I$ is smaller than $J$ for the lexicographic 
order on $\PC(S_{n-1})$ induced by the order $t < s_1 < \dots < s_{n-1}$ 
on $S_{n-1}$. 
\end{itemize}
\end{quotation}
It follows immediately from $(*)$ that 
$$\Res_{S_{n-1}}^{S_n} x_J^{S_n} \in \a_J x_J^{S_{n-1}} + 
\sum_{I \inferieur J} \QM x_I^{S_{n-1}}$$
with $\a_J > 0$ (for every $J \in \PC(S_{n-1})$). The proof of the proposition 
is complete.
\end{proof}

\bigskip

\begin{coro}\label{surjectif D}
The image of the map 
$\Res_{S_{n-1}'}^{S_n'} : \Sigma(W_n') \to \Sigma(W_{n-1}')$ 
is equal to $\Sigma(W_{n-1}')^{\s_{n-1}}$.
\end{coro}

\begin{proof}
This follows from Proposition \ref{restriction BD} (c) and (e) 
and from Proposition \ref{surjectif B}.
\end{proof}

\remark{sigma n}
If $n$ is odd, then $\s_n=\s_0$, the automorphism of $\Sigma(W_n')$ 
induced by conjugation by the longest element of $W_n'$.\finl

\bigskip

\soussection{Self-opposed subsets} 
A subset $K$ of $S$ is called {\it self-opposed} if, 
for every $w \in W$ such that $\lexp{w}{K} \subset S$, we have 
$\lexp{w}{K}=K$. 

In this subsection, we fix a self-opposed 
subset $K$ of $S$. 
If $s \in S\setminus K$, we set $w_{K,s} = w_{K \cup \{s\}} w_K$ 
(here, if $I$ is a subset of $S$, $w_I$ denotes the longest 
element of $W_I$). Then, since $K$ is self-opposed, we have $w_{K,s} \in W(K)$. 
Now, if $I$ is a subset of $S$ containing $K$, we set 
$I(K)=\{w_{K,s}~|~s \in I \setminus K\}$. Then (see 
for instance \cite[Remark 2.3.5]{geck pfeiffer})
\equat\label{auto oppose}
(W(K),S(K))\text{\it~ is a finite Coxeter group.}
\endequat

\medskip

\begin{quotation}
\noindent{\sc Notation - } 
We denote by $X_I^{(K)}$, $x_I^{(K)}$, $\L_{(K)}$ and $\lamb_{(K)}$ 
the objects defined like $X_I$, $x_I$, $\L$ or $\lamb$ but inside $W(K)$.  
\end{quotation}

\medskip

Let $\psi_K : \Sigma(W) \to \Sigma(W(K))$ be the linear map such that 
$$\psi_K(x_I) = \begin{cases}
                x_{I(K)}^{(K)} & \text{if $K \subset I$,} \\
                0 & \text{otherwise,}
                \end{cases}$$
for every subset $I$ of $S$.

\medskip

\begin{prop}\label{goetz 1}
  Assume $K$ is self-opposed in $S$.
  Let $I, J, L \subseteq S$ be  such that $K \subseteq I, J, L$.  Then
  $X^{(K)}_{I(K) J(K) L(K)} = X_{IJL}$.
\end{prop}

\begin{proof}
  First note that, if $l_{(K)}$ is the length function of $W(K)$
  with respect to $S(K)$ then, for any $s \in S \setminus K$, we
  have that $l(ws) > l(w)$ if and only if $l_{(K)}(w w_{K,s}) >
  l_{(K)}(w)$ (see \cite[Theorem 5.9]{Lusztig1976}).  It follows
  that $X^{(K)}_{J(K)} = X_J \cap W(K)$ for every subset $J$ of $S$
  containing $K$.  Moreover, $W_{J(K)} = W_J \cap W(K)$.

  Let $d \in X_{IJL}$ for some $L \subseteq S$ containing $K$.  Then
  $K \subseteq L \subseteq I^d$ implies $d \in W(K)$.  Also $W_I^d
  \cap W_J = W_L$ implies $W_{I(K)}^d \cap W_{J(K)} = (W_I \cap
  W(K))^d \cap W_J \cap W(K) = W_L \cap W(K) = W_{L(K)}$, whence
  $X_{IJL} \subseteq X^{(K)}_{I(K) J(K) L(K)}$ for all $K \subseteq L
  \subseteq S$.

  Equality follows from that fact that $X^{(K)}_{I(K) J(K)} = X_{IJ}
  \cap W(K)$ is both the disjoint union of the sets $X^{(K)}_{I(K)
    J(K) L(K)}$ with $K \subseteq L \subseteq S$ and the disjoint
  union of the sets $X_{IJL}$ with $K \subseteq L \subseteq S$.
\end{proof}

\begin{theo}\label{goetz}
If $K$ is a self-opposed subset of $S$, then $\psi_K$ is a surjective 
homomorphism of algebras.
\end{theo}

\begin{proof}
The surjectivity of $\psi_K$ is clear from the definition. Also, 
$\psi_K(1)=\psi_K(x_S)=x_{S(K)}=1$. 
Let us now prove that $\psi_K$ is respects the multiplication. 
Let $I$, $J$ be two subsets of $S$. We want to prove that 
$$\psi_K(x_Ix_J)=\psi_K(x_I)\psi_K(x_J).\leqno{(*)}$$

Assume first that $I$ (or $J$) does not contain $K$. Then 
$\psi_K(x_I)\psi_K(x_J)=0$. Let $d \in X_{IJ}$. 
If $K$ is contained in $\lexp{d^{-1}}{I} \cap J$, then 
$K$ is contained in $J$ or $\lexp{d}{K}$ is contained in $I \subset S$, 
so $K$ is contained in $I$ or in $J$, which is impossible. So 
$\psi_K(x_Ix_J)=0$.

Assume now that both $I$ and $J$ contain $K$. 
Then
$$\psi_K(x_Ix_J)= \sum_{K \subset L \subset S} |X_{IJL}|~x_{L(K)}^{(K)}.$$
In this case, $(*)$ follows from Proposition \ref{goetz 1}. 
\end{proof}

\medskip

\example{Bn}
Assume here that $W$ is of type $B_n$ and keep the notation 
of the proof of Proposition \ref{surjectif B}. Then $\{t\}$ 
is a self-opposed subset of $S$ and $W(\{t\})$ is of type $B_{n-1}$. 
So Theorem \ref{goetz} gives another surjective morphism between 
the Solomon algebra of type $B_n$ and the Solomon algebra of type 
$B_{n-1}$. This homomorphism does not coincide with 
the one constructed in Proposition \ref{surjectif B}.\finl

\medskip

\example{E7} 
Assume here that $(W,S)$ is of type $E_7$ and assume 
that $S=\{s_i~|~1 \le i \le 7\}$ 
is numbered as in \cite[Planche VI]{bourbaki}. In other words, the Dynkin 
diagram of $W$ is:
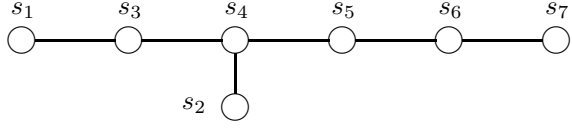
\begin{center}
\begin{picture}(240,45)
\put( 40, 30){\circle{10}}
\put( 45, 30){\line(1,0){30}}
\put( 80, 30){\circle{10}}
\put( 85, 30){\line(1,0){30}}
\put(120, 30){\circle{10}}
\put(125, 30){\line(1,0){30}}
\put(160, 30){\circle{10}}
\put(165, 30){\line(1,0){30}}
\put(200, 30){\circle{10}}
\put(205, 30){\line(1,0){30}}
\put(240, 30){\circle{10}}
\put(120, 5){\circle{10}}
\put(120, 25){\line(0,-1){15}}
\put( 36, 40){$s_1$}
\put( 76, 40){$s_3$}
\put(116, 40){$s_4$}
\put(156, 40){$s_5$}
\put(196, 40){$s_6$}
\put(236, 40){$s_7$}
\put(100, 4){$s_2$}
\end{picture}
\end{center}
Let $K=\{s_2,s_5,s_7\}$. Then $K$ is self-opposed and $W(K)$ is of type 
$F_4$. So Theorem \ref{goetz} realizes the Solomon algebra of type 
$F_4$ as a quotient of the Solomon algebra of type $E_7$.\finl

\bigskip

If $I$ is a subset of $S(K)$, we denote by $\varpi_K(I)$ 
the unique subset $A$ of $S$ containing $K$ such that $A(K)=I$. 
Then the map $\varpi_K : \PC(S(K)) \to \PC(S)$ induces 
a map $\varpit_K : \L_{(K)} \to \L$ (indeed, by the definition 
of $W(K)$, if $I$ and $J$ are two subsets of $S(K)$ and if $w \in W(K)$ 
is such that $\lexp{w}{I}=J$, then $\lexp{w}{\varpi_K(I)}=\varpi_K(J)$). 
Then, if $I \subset S(K)$, we have, by \ref{tau bis},
\equat\label{varpi}
\t_{\lamb(\varpi_K(I))}=\t_{\lamb_{(K)}(I)}^{(K)}  \circ \psi_K.
\endequat

We close this subsection by showing that the morphisms $\Res_L$ and 
$\psi_K$ are compatible. More precisely, let $L$ be a subset 
of $S$ containing $K$. Then $K$ is self-opposed for $W_L$ and 
$W_L(K)$ is the parabolic subgroup of $W(K)$ generated by $L(K)$. 
Let $\psi_K^L : \Sigma(W_L) \to \Sigma(W_L(K))$ be the morphism 
defined like $\psi_K$ but inside $W_L$. Then 
the diagram
\equat\label{commute phi psi}
\diagram
\Sigma(W) \rrto^{\DS{\psi_K}} \ddto_{\DS{\Res_L}} && 
\Sigma(W(K)) \ddto^{\DS{\Res_{L(K)}}} \\
&& \\
\Sigma(W_L) \rrto^{\DS{\psi_K^L}} && \Sigma(W_L(K)) \\
\enddiagram
\endequat
is commutative. Indeed, if $I$ is a subset of $L$ containing $K$, 
we have $\psi_K(x_I)=x_{I(K)}=x_{L(K)} \psi_K^L(x_I^L)$. In other words, 
$\psi_K(x_L x) = x_{L(K)} \psi_K^L(x)$ for every $x \in \Sigma(W_L)$. 
So the commutativity of \ref{commute phi psi} 
follows from Proposition \ref{bbht} (a) and 
from routine computations.

\bigskip

\section{Loewy length of $\Sigma(W)$}

\medskip

The {\it Loewy length} of a finite dimensional algebra $A$ is the smallest 
natural number $k \ge 1$ such that $(\Rad A)^k = 0$. 
We denote by $\loewy(W)$ the Loewy length of $\Sigma(W)$. If 
$\sigma$ is an automorphism of $W$ such that $\s(S)=S$, we denote by 
$\loewy(W,\s)$ the Loewy length of $\Sigma(W)^\s$. By Corollary 
\ref{radical fixe}, we have
\equat\label{inegalite loewy}
\loewy(W,\s) \le \loewy(W).
\endequat
By Solomon's Theorem (e), $\loewy(W)$ is the smallest natural number 
$k \ge 1$ such that $(\Ker \th)^k=0$. 

\bigskip

\soussection{Upper bound} 
Let us start with an easy observation (recall that $\s_0$ denotes 
the automorphism of $W$ induced by conjugacy by $w_0$):

\begin{lem}\label{diminution}
Let $k \ge 0$. Then:
\begin{itemize}
\itemth{a} $(\Ker \th) . \Sigma_k(W) \subset \Sigma_{k-1}(W)$. 

\itemth{b} $(\Ker \th)^{\s_0} . \Sigma_k(W) \subset \Sigma_{k-2}(W)$. 
\end{itemize}
\end{lem}

\begin{proof}
Let $J \in \PC(S)$ be such that $|J|\le k$ and let $x \in \Ker \th$. 
Then $\t_{\lamb(J)}(x)=0$. By \ref{tau bis}, we then have 
$x x_J \in \Sigma_{k-1}(W)$, whence (a). If moreover $x \in (\Ker \th)^{\s_0}$, 
then $x x_J' \in \Sigma_{k-2}(W)$ by Lemma \ref{w0 central}. 
This shows (b).
\end{proof}

\bigskip

\noindent{\sc Remark - } It is not true in general that 
$\Sigma_k(W).(\Ker \th)   \subset \Sigma_{k-1}(W)$.\finl

\begin{coro}\label{borne loewy}
We have:
\begin{itemize}
\itemth{a} $\loewy(W) \le |S|$.

\itemth{b} $\loewy(W,\s_0) \le \DS{\frac{|S|+1}{2}}$.
\end{itemize}
\end{coro}

\begin{proof}
(a) We have $\Ker \th \subset \Sigma_{|S|-1}(W)$ and 
$\Ker \th \cap \Sigma_0(W) = 0$ (see Solomon's Theorem (d)). So, by 
Lemma \ref{diminution} (a), we have $(\Ker \th)^{|S|} = 0$. 

\medskip

(b) By Lemma \ref{diminution} (b), we have 
$(\Ker \th)^{\s_0} \subset \Sigma_{|S|-2}(W)$ and $\bigl((\Ker \th)^{\s_0}\bigr)^r 
\subset \Sigma_{|S|-2r}(W)$ for every $r \ge 0$. This shows (b).
\end{proof}

\medskip

\example{loewy A} 
It is a classical result \cite[Corollary 3.5]{atkinson} that, 
if $W$ is of type $A_n$, then $\loewy(W)=n$. In this case, we also have 
$\loewy(W,\s_0)=\DS{\Bigl\lceil\frac{n}{2}\Bigr\rceil}$. Indeed, let $l=\loewy(W,\s_0)$. 
By Corollary \ref{borne loewy} (b), we have 
$l \le \DS{\Bigl\lceil\frac{n}{2}\Bigr\rceil}$. 
On the other hand, let $a = x_{\{s_1,\dots,s_{n-1}\}}-x_{\{s_2,\dots,s_n\}}$, 
where $S=\{s_1,s_2,\dots,s_n\}$ is numbered 
such that $(s_is_{i+1})^3=1$ for every 
$i \in \{1,2,\dots,n-1\}$. Then, by \cite[Proof of Corollary 3.5]{atkinson}, 
we have $a \in \Rad \Sigma(W)$ and $a^{n-1} \neq 0$. In particular, 
$(a^2)^{\SS{\bigl[\frac{n-1}{2}\bigr]}} \neq 0$. 
But, $\s_0(a)=-a$, so $\s_0(a^2)=a^2$. Therefore, by Corollary \ref{radical fixe}, 
we have $a^2 \in \Rad\bigl(\Sigma(W)^{\s_0}\bigr)$. So 
$l \ge \DS{\Bigl\lceil\frac{n}{2}\Bigr\rceil}$, as desired.\finl

\bigskip

\soussection{Type B} 
We keep the notation of subsection \ref{section BD}. 
The aim of this subsection is to prove the next 
proposition:

\begin{prop}\label{loewy B}
If $n \ge 1$, then $\loewy(W_n)=\DS{\Bigl\lceil\frac{n}{2}\Bigr\rceil}$.
\end{prop}

\begin{proof}
By Corollary \ref{borne loewy} (b), we have 
$\loewy(W_n)\le\DS{\Bigl\lceil\frac{n}{2}\Bigr\rceil}$. Now, let 
$r=\DS{\Bigl[\frac{n-1}{2}\Bigr]}$. It is sufficient to find 
$a_1$,\dots, $a_r \in \Rad \Sigma(W_n)$ such that $a_r \dots a_1 \neq 0$. 

If $1 \le i \le j \le n-1$, we set $[i,j]=\{s_i,s_{i+1},\dots,s_j\}$. 
If $1 \le i \le r$, we set 
$$a_i = x_{[2i-1,n-2]}-x_{[2i,n-1]}.$$
Then $a_i \in \Rad \Sigma(W_n)$. We shall show that $a_r \dots a_1 \neq 0$. 
If $1 \le i \le r$, 
we set 
$$\t_i = \sum_{j=0}^{2i-1} (-1)^j \matrice{2i-1 \\ j} x_{[j+1,n-2i+j]}.$$
We will show by induction on 
$i$ that 
$$a_i \dots a_1 \in \QM^\times \t_i + \Sigma_{n-2i-1}(W_n).\leqno{(P_i)}$$
Note that, if $(P_r)$ is proved, then 
the proposition is complete. Now, $(P_1)$ holds since $a_1=\t_1$. 
So, let $i \in \{2,3,\dots,r\}$ and assume that $(P_{i-1})$ holds. 
By Lemma \ref{diminution} (b), 
there exists three elements $\a$, $\b$ and $\g$ of $\QM$ such that 
\begin{multline*}
a_i x_{[j+1,n-2(i-1)+j]} \\ \in \a x_{[j+1,n-2i+j]} + \b 
x_{[j+2,n-2i+j+1]} +\g x_{[j+3,n-2i+j+2]} \\+ \Sigma_{n-2i-1}(W)
\end{multline*}
for every $j \in \{1,2,\dots 2i-1\}$. The fact that $\a$, $\b$ and $\g$ 
do not depend on $j$ follows from the fact that there exists 
$w \in X_{[j+1,n-2(i-1)+j],[j'+1,n-2(i-1)+j']}$ such that 
$\lexp{w}{[j'+1,n-2(i-1)+j']}=[j+1,n-2(i-1)+j]$. In particular, 
we have 
$x_{[j+1,n-2(i-1)+j]}w=x_{[j'+1,n-2(i-1)+j']}$. 
Since $a_i \in \Ker \th$, we have $\a+\b+\g=0$. Also, 
$$x_{[j+1,n-2(i-1)+j]}w_0w_{[j+1,n-2(i-1)+j]}=x_{[j+1,n-2(i-1)+j]}.$$
Therefore, $\a=\g$. In other words, 
\begin{multline*}
a_i x_{[j+1,n-2(i-1)+j]} \\
\in \a (x_{[j+1,n-2i+j]} -2 
x_{[j+2,n-2i+j+1]} + x_{[j+3,n-2i+j+2]}) \\
+ \Sigma_{n-2i-1}(W).
\end{multline*}
Hence, by the induction hypothesis, by Lemma \ref{diminution} (b) 
and by usual properties of binomial coefficients, we have 
$$a_i \dots a_1 \in \QM^\times \a\t_i + \Sigma_{n-2i-1}(W_n).$$
So it remains to show that 
$$\a \neq 0.$$
For this, consider the case where $j=2i-3$ and write
$$x_{[2i-1,n-2]} x_{[2i-2,n-1]} = a x_{[2i-2,n-3]}+bx_{[2i-1,n-2]}+cx_{[2i,n-1]}$$
$$x_{[2i,n-1]} x_{[2i-2,n-1]}= d x_{[2i-2,n-3]}+ex_{[2i-1,n-2]}+fx_{[2i,n-1]}
\leqno{\text{and}}$$
with $a$, $b$, $c$, $d$, $e$ and $f$ in $\QM$. 
We then have $b-e = -2\a$. Since $b \not= 0$, it is sufficient to 
show that $e=0$. In other words, we need to prove the following lemma:

\begin{quotation}
\begin{lem}\label{fin}
If $d \in X_{[2i,n-1],[2i-2,n-1]}$, then 
$\lexp{d^{-1}}{[2i,n-1]}\not=[2i-1,n-2]$.
\end{lem}

\begin{proof}[Proof of Lemma \ref{fin}]
We identify $W_n$ with the group of permutations $\s$ 
of $E=\{\pm 1, \pm 2, \dots, \pm n\}$ such that $\s(-k)=-\s(k)$ 
for every $k \in E$ ($t$ corresponds to the transposition $(-1,1)$ 
while $s_k$ corresponds to $(k,k+1)(-k,-k-1)$). If 
$d \in X_{[2i-2,n-1]}$, then $d$ is increasing on 
$\{2i-2,2i-1,\dots,n-2,n-1\}$. If moreover 
$\lexp{d^{-1}}{[2i,n-1]}=[2i-1,n-2]$ (in other words, if 
$\lexp{d}{[2i-1,n-2]}=[2i,n-1]$), then 
$d(\{2i-1,2i,\dots,n-1\}) \subset \{\pm 2i, \pm (2i+1), \dots \pm n\}$ 
and $d$ has constant sign on $\{2i-1,2i,\dots,n-1\}$.
Two cases may occur:

$\bullet$ If $d(2i-1) > 0$, then, since $d$ has constant sign and is increasing 
on $\{2i-1,2i,\dots,n-1\}$, we have $d(n-1)=n$. But this is impossible 
since $d(n) > d(n-1)$.

$\bullet$ If $d(2i-1) < 0$, then, by the same argument, we have $d(2i-1) = -n$. 
But this is again impossible since $d(2i-2) < d(2i-1)$.
\end{proof}
\end{quotation}

The proof of $(P_i)$ and of the proposition is now complete.
\end{proof}

\bigskip

\remark{b tau}
Assume here that $W$ is of type $B_{2r+1}$, $r \ge 1$ and let $\t_r$ denote 
the element of $\Sigma(W)$ defined in the proof of Proposition \ref{lower bound}. 
Computations using {\tt CHEVIE} show that the following question has a positive 
answer for $m \in \{1,2,3\}$:

\begin{quotation}
\noindent{\bf Question:} {\it Is it true that $(\Ker \th)^r = \QM \t_r$?\finl}
\end{quotation}

\bigskip

\bigskip

\soussection{Type D}
The following result is an easy consequence of Proposition 
\ref{loewy B} (and its proof) and of the existence of the homomorphism 
of algebras $\Res_n$.

\begin{coro}\label{loewy D fixe}
If $n \ge 1$, then $\loewy(W_n',\s_n)=\DS{\Bigl\lceil\frac{n}{2}\Bigr\rceil}$. 
\end{coro}

\begin{proof}
Let $l$ be the Loewy length of $\Sigma(W_n')^{\s_n}$. 
Keep the notation of the proof of Proposition \ref{loewy B}. 
Let $a_i'=\Res_n a_i$. Then, by Proposition \ref{restriction BD} 
(b) and (e), we have $a_i' \in \Rad \Sigma(W_n')^{\s_n}$ and 
$a_r\dots a_1 = \a \Res_n \t_r$, where $\a \not= 0$. But, 
it is clear from the definition of $\Res_n$ that 
$\Res_n \t_r \neq 0$. So $a_r' \dots a_1' \neq 0$. So 
$l \ge r+1$. 

The fact that $l \le r+1$ follows from Propositions \ref{loewy B} and 
\ref{restriction BD} (e).
\end{proof}

\begin{coro}\label{loewy D}
Let $n \ge 3$. Then:
\begin{itemize}
\itemth{a} If $n$ is even, then $\loewy(W_n') = \DS{\Bigl\lceil\frac{n}{2}\Bigr\rceil}$. 

\itemth{b} If $n$ is odd, then $\loewy(W_n') \ge \DS{\frac{n+3}{2}}$. 
\end{itemize}
\end{coro}

\begin{proof}
By Proposition \ref{restriction BD} (e), 
$$\Res_n(\Rad \Sigma(W_n)) = \Rad \Sigma(W_n')^{\s_n} = 
(\Rad \Sigma(W_n'))^{\s_n}.\leqno{(*)}$$
So (a) follows from $(*)$, from Corollary \ref{borne loewy} 
and from Corollary \ref{loewy D fixe}.

Let us now prove (b). Write $n=2r+1$ and keep the notation of the proof 
of Corollary \ref{loewy D fixe}. 
Let $a=x_{\{s_1',s_2,\dots,s_{2r}\}}-x_{\{s_1,s_2,\dots,s_{2r}\}}$. 
An easy computation shows that 
$$a_1'  x_{\{s_1,s_2,\dots,s_{2r}\}} \in  a_1'
+\Sigma_{\PC^\#(\{s_1,s_2,\dots,s_{2r}\})}(W).$$
But, by the equalities (3) and (4) of the proof of Proposition 
\ref{lower bound} and by Lemma \ref{diminution} (b), we have 
$\s_r \dots \s_2 \Sigma_{\PC^\#(\{s_1,s_2,\dots,s_{2r}\})}(W)=0$. 
Therefore, 
$$\s_r\dots\s_1 x_{\{s_1,s_2,\dots,s_{2r}\}} = \s_r\dots\s_1.$$ 
Since $x_{\{s_1',s_2,\dots,s_{2r}\}}=x_{\{s_1,s_2,\dots,s_{2r}\}} d$, 
where $d=w_{[1,r]}w_0$, we get that 
$$\s_r\dots\s_1 a = \s_r\dots\s_2\s_1(1-d).$$
Therefore, $\s_r \dots \s_1 a \not= 0$ (indeed, the coefficient of 
$x_{s_1'}$ is non-zero). But, $a \in \Ker \th$ because $|S|$ is odd. 
So $\loewy(W_{2r+1}') \ge r+2$, as desired.
\end{proof}

\bigskip

\soussection{Lower bound} 
The aim of this subsection is to prove the 
following result:

\begin{prop}\label{lower bound}
If $W$ is irreducible, then $\loewy(W) \ge \loewy(W,\s_0) \ge 
\DS{\Bigl\lceil\frac{|S|}{2}\Bigr\rceil}$.
\end{prop}

\begin{proof}
By \ref{inegalite loewy}, we only need to prove the second inequality. 
The proof of this proposition will proceed by a case-by-case analysis. 
First, the exceptional groups can be treated by using {\tt CHEVIE} 
(see the Table given at the end of this paper). 
If $W$ is of type $A$, then this follows from Example \ref{loewy A}. 
If $|S|=2$, then there is nothing to prove. 
If $W$ is of type $B$, this follows from Proposition \ref{loewy B}. 
If $W$ is of type $D$, this follows from Corollary \ref{loewy D}. 
The proof is complete. 
\end{proof}

The next result follows from Corollary \ref{borne loewy} 
and Proposition \ref{lower bound} 

\begin{coro}\label{loewy}
If $W$ is irreducible, then 
$\loewy(W,\s_0)=\DS{\Bigl\lceil\frac{|S|}{2}\Bigr\rceil}$. 
\end{coro}

\begin{coro}\label{loewy central}
If $W$ is irreducible and $w_0$ is central in $W$, then 
$\loewy(W)=\DS{\Bigl\lceil\frac{|S|}{2}\Bigr\rceil}$. 
\end{coro}

\bigskip

\soussection{Conclusion} The next table gives the known Loewy lengths 
of the algebras $\Sigma(W)^\s$ for $W$ irreducible and 
$\s$ is a length-preserving automorphism of $W$. 
$$
\begin{array}{|c|c||c|c|l|}
\hline
\text{Type of }W & o(\s) & |\L/\s| & \loewy(W,\s) & d_0, d_1,d_2,\dots \verticalhb\\
\hline
\hline
A_n \verticalh & 1 & p(n) & n & \\
  \verticalhbb   & 2 & p(n) & \Bigl\lceil\DS{\frac{n}{2}}\Bigr\rceil & \\
\hline
B_n \verticalhbb & 1 & \DS{\sum_{r=0}^n p(r)} & \Bigl\lceil\DS{\frac{n}{2}}\Bigr\rceil & \\
\hline
D_{2n} \verticalh & 1 & \DS{p(n) + p(2n) + \sum_{r=0}^{2n-2} p(r)} & n & \\
(n \ge 2) & 2 & \DS{p(2n) + \sum_{r=0}^{2n-2} p(r)} & n & \verticalhb \\
\hline
D_{2n+1} \verticalh & 1 & \DS{p(2n) + \sum_{r=0}^{2n-2} p(r)} & \ge n+2~^{(*)} & \\
(n \ge 2) & 2 & \DS{p(2n) + \sum_{r=0}^{2n-2} p(r)} & n+1 & \verticalhb \\
\hline
D_4 \verticalh & 1 & 11 & 2 & 16, 5  \\
    \verticalh & 2 & 9 & 2 & 12, 3 \\
    \verticalh & 3 & 7 & 2 & 8, 1 \\
\hline
E_6 & 1 & 17 & 5 & 64,47,28,12,3 \verticalh \\
    & 2 & 17 & 3 & 40, 23, 5  \verticalhb \\
\hline
E_7 & 1 & 32 & 4 & 128,96,34,2 \verticalhb \\
\hline
E_8 & 1 & 41 & 4 & 256,215,106,14 \verticalhb \\
\hline
F_4 & 1 & 12 & 2 & 16, 4 \verticalh \\
    & 2 & 8  & 2 & 10, 2 \verticalhb \\
\hline
H_3 & 1 & 6  & 2 & 8,2 \verticalhb \\
\hline
H_4 & 1 & 10 & 2 & 16,6 \verticalhb \\
\hline
I_2(2m)   & 1 & 4 & 1 & 4  \verticalh \\
          & 2 & 3 & 1 & 3 \verticalhb \\
\hline
I_2(2m+1) & 1 & 3 & 2 & 4,1 \verticalh \\
          & 2 & 3 & 1 & 3 \verticalhb \\
\hline
\end{array}
$$
The exceptional groups are obtained by using {\tt CHEVIE}. 
Type A is mainly due to Atkinson \cite[Corollary 3.5]{atkinson} 
(see Example \ref{loewy A}). Types $B$ and $D$ are done in this paper 
(except for the type $D_{2n+1}$). Dihedral groups 
are easy. It must be noticed that the inequality $(*)$ is an equality 
for $n=2$ or $3$. We suspect it is always an equality. 

In this table, 
$d_i$ denotes the dimension of $\bigl(\Rad(\Sigma(W)^\s)\bigr)^i$. 
These numbers are not given for infinite series of type $A$, $B$ or $D$. 
Note that $d_0=\dim \Sigma(W)^\s$ and that 
$|\L/\s|=d_0-d_1$. We denote by $p(n)$ the number of 
partitions of $n$. We denote by $o(\s)$ the order of $\s$: 
it characterizes the conjugacy class of $\s$ in the group 
of automorphism of $W$ stabilizing $S$. 

%
%

\end{document}